\documentclass[reqno]{amsart}

\numberwithin{equation}{section}

\newtheorem{Thm}{Theorem}
\newtheorem{Lem}[Thm]{Lemma}

\newtheorem{Prop}{Proposition}

\theoremstyle{definition}
\newtheorem{Def}[Thm]{Definition}

\theoremstyle{remark}
\newtheorem*{Rem}{Remark}
\theoremstyle{remark}

\newtheorem*{Notn}{Notation}

\def\demo#1{\vskip-\lastskip
\vskip12pt \noindent{\it#1.\/} }
\def\enddemo{\vskip12pt}

\newcommand{\e}{\varepsilon}
\newcommand{\G}{\Gamma}
\newcommand{\Z}{\mathbf{Z}}

\newcommand{\ds}{\displaystyle}
\newcommand{\bigoh}{\mathrm{O}}
\newcommand{\il}{\int_0^\infty}

\input epsf

\def\hint#1{\int\limits_{\epsffile{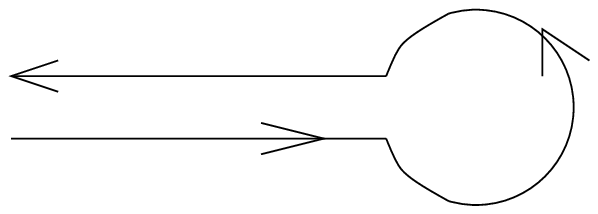}}\!\!\!#1}

\begin{document}
\title[Difference Differential Equations]{A Pair of Difference
Differential Equations of Euler-Cauchy Type}

\author{David M. Bradley}
\address{Department of Mathematics and Statistics\\
         University of Maine\\
         5752 Neville Hall\\
         Orono, Maine 04469--5752\\
         U.S.A.}
\email{dbradley@member.ams.org, bradley@math.umaine.edu}
\thanks{Research supported by the University
of Maine summer faculty research fund.}
\thanks{Published in \textit{Transactions of the American Mathematical
Society}, \textbf{355} (2003), no.~12, 4985--5002.
http://www.ams.org/journal-getitem?pii=S0002-9947-03-03223-9}


\subjclass{Primary: 34K06; Secondary: 34K12, 34K25}

\keywords{Difference differential equations, integral transforms,
adjoint relation, Dickman-de Bruijn function, sieves.}

\begin{abstract}
We study two classes of linear difference differential equations
analogous to Euler-Cauchy ordinary differential equations, but in
which multiple arguments are shifted forward or backward by fixed
amounts.  Special cases of these equations have arisen in diverse
branches of number theory and combinatorics.  They are also of use
in linear control theory.  Here, we study these equations in a
general setting.  Building on previous work going back to de
Bruijn, we show how adjoint equations arise naturally in the
problem of uniqueness of solutions.  Exploiting the adjoint
relationship in a new way leads to a significant strengthening of
previous uniqueness results.  Specifically, we prove here
(Theorem~\ref{thm:exp}) that the general Euler-Cauchy difference
differential equation with advanced arguments has a unique
solution (up to a multiplicative constant) in the class of
functions bounded by an exponential function on the positive real
line.  For the closely related class of equations with retarded
arguments, we focus on a corresponding class of solutions,
locating and classifying the points of discontinuity. We also
provide an explicit asymptotic expansion at infinity.
\end{abstract}

\maketitle

\section{Introduction}
We study two classes of linear difference differential equations
analogous to Euler-Cauchy ordinary differential equations, but in
which multiple arguments are shifted forward or backward by fixed
amounts. Special cases of these equations have been studied quite
extensively in the context of certain problems in number theory
and combinatorics---see \S\ref{sect:back} for a brief overview.
Here, we take a broader viewpoint.  As in~\cite{BraDi}, we focus
primarily on the advanced-argument linear difference differential
equation
\begin{equation}
   uq'(u)=\sum_{j=0}^m \alpha_j q(u+v_j),\qquad u>0,
   \label{+ddereal}
\end{equation}
in which the argument $u$ is incremented by the non-negative real
numbers $0=v_0<v_1<\cdots<v_m$, and the coefficients
$\alpha_0,\alpha_1,\ldots,\alpha_m$ are arbitrary complex numbers.
Here, however, we study~(\ref{+ddereal}) in conjunction with the
related delay differential equation
\begin{equation}
   (up(u))'=-\sum_{j=0}^m \alpha_j p(u-v_j).
   \label{-ddereal}
\end{equation}
It turns out that by judiciously pairing equations from each
class, one can infer properties of solutions to one equation from
properties of a solution to the other.  We exploit this
relationship to significantly strengthen previous uniqueness
results concerning equation~(\ref{+ddereal}). Our main uniqueness
result is that~(\ref{+ddereal}) has a unique solution (up to a
multiplicative constant) in the class of exponentially bounded
functions.

To exploit the relationship between the two
equations~(\ref{+ddereal}) and~(\ref{-ddereal}), we need to amass
a certain amount of information concerning a solution to the
latter.  In this direction we establish new results, including the
nature and location of discontinuities, and an explicit asymptotic
expansion.

The organization of the paper is as follows.  In
\S\ref{sect:back}, a brief historical overview is provided.  We
discuss uniqueness of solutions to~(\ref{+ddereal}) in the class
of polynomially bounded functions in \S\ref{sect:poly}.
The next section focuses on a special case of~(\ref{+ddereal})
arising in sieve theory and the study of cycle lengths of a random
permutation. The particular properties of the parameters in this
case are exploited to show that any solution to the underlying
equation which is not polynomially bounded must be very wild
indeed.  The techniques of \S\ref{sect:sieveaux} prepare the
ground for our exponential uniqueness theorem.  We motivate the
introduction of the adjoint equation in \S\ref{sect:adj}.
Information concerning the asymptotic behavior of solutions to the
adjoint equation is needed, and the relevant results are presented
in \S\ref{sect:pasym}.  Details of the proof of our exponential
uniqueness theorem are relegated to the penultimate section.

\begin{Notn} The
Bachmann-Landau $\bigoh$-notation and the standard notion of
asymptotic equivalence are briefly recalled in the forms that we
wish to use them here. If $g$ is a positive function of a positive
real variable, the symbol $\bigoh(g(u))$, $0<u\to\infty$ denotes
an unspecified function $f$ for which there exist positive real
numbers $u_0$ and $B$ such that $|f(u)|\le Bg(u)$ for all real $u>
u_0$.  Two functions $f$ and $g$ are asymptotically equivalent at
infinity, written $f(u)\sim g(u)$, $0<u\to\infty$, if for every
$\e>0$, there exists a positive real number $u_0$ such that
$|f(u)-g(u)|<\e |g(u)|$ whenever $u>u_0$.
\end{Notn}

\section{Historical Overview}
\label{sect:back} Ever since Dickman~\cite{Di} showed that the
asymptotic density of the positive integers $n\le x$ having no
prime factors exceeding $x^{1/u}$ is given by the continuous
solution $\rho(u)$ to the difference differential equation
\begin{alignat}{2}
   u\rho'(u) &=-\rho(u-1), \qquad & u>1,\notag
   \intertext{with boundary condition}
   \rho(u)   &=1,  \qquad & 0<u\le 1,\notag
\end{alignat}
it has been known that difference differential equations arise in
the study of certain problems in analytic number theory. We now
know of at least four types of number theoretical problems in
which difference differential equations arise.  These may be
loosely classified according to the context of the original
problem. The first type, ``psixyology''~\cite{Moree} encompasses
those problems associated with determining the probability
distribution of the prime factors of an
integer~\cite{Bi,dB2,dB3,Bu,Di,Hi}. The second type, ``sieve
theory''~\cite{HRbook} is concerned with the problem of estimating
the size of a finite set of integers after certain residue classes
have been
eliminated~\cite{AO,dB1,CG,DHR1,DHR1II,DHR2,DHR3,DHR4,Gr,Iw}. The
third type deals with the problems relating to the determination
of the cycle decomposition of a random
permutation~\cite{AB,Bal,Gol,Gonch,KTP,PW,SL}, and is closely
related to the first~\cite{AMSNotices}.  The fourth type focuses
on estimating incomplete sums of multiplicative
functions~\cite{Al,dBvL,He}.

Despite the extensive literature on the subject, apart from the
treatment of a few isolated
cases~\cite{DBrad0,DBrad1,DHR5,DHR6,DHR7,Hi2,HiTe,tRi,Wh1,Wh2} the
main focus has been on the application to problems in number
theory and combinatorics as opposed to the study of the underlying
difference differential equations as a problem in its own right.
Although a systematic study of equation~(\ref{+ddereal}) in full
generality was initiated in~\cite{BraDi}, there is still
additional work that remains to be done with regard to providing a
comprehensive treatment. The present paper can be viewed as a
further contribution in this respect, in that we explore more
deeply the question of uniqueness, and reveal more fully the role
of the adjoint equation~(\ref{-ddereal}).  It is hoped that some
of the analytical techniques used herein---and more generally our
ongoing program of study---will be of interest to researchers in
number theory and the wider differential equations community as
well.

%
%
%

\section{Polynomially Bounded Solutions}
\label{sect:poly}

As in~\cite{BraDi}, we fix a non-negative integer $m$, real
numbers $0=v_0<v_1<\dots < v_m$, complex numbers $\alpha_0,
\alpha_1,\dots, \alpha_m$, and put
$\beta:=\alpha_0+\alpha_1+\cdots+\alpha_m$.  Without loss of
generality, we assume $\alpha_j\ne 0$ for $1\le j\le m$.  Define a
function $q^*$ on the positive real half-line $(0..\infty)$ as
follows.
\begin{Def}\label{def:q*}   For $u>0$, let
\begin{subequations}\label{q*}
\begin{align}
   q^*(u) &:= \ds\frac{1}{\G(-\beta)}\il x^{-\beta-1}
   \exp\bigg\{-ux+\sum_{j=1}^m\alpha_j\int_0^x
       \frac{1-e^{-v_j t}}{t}\,dt\bigg\}\,dx, \label{qLaplace}\\
          &\text{if $\beta$ is a negative integer, and}\notag\\
    q^*(u) &:= \ds\frac{\G(\beta+1)}{2\pi i}\hint{z^{-\beta-1}}
      \exp\bigg\{uz-\sum_{j=1}^m\alpha_j
      \int_0^z\dfrac{e^{v_jt}-1}{t}\,dt\bigg\}\,dz, \label{qHankel}\\
          &\text{otherwise.}\notag
\end{align}
\end{subequations}
\end{Def}
The contour in~(\ref{qHankel}) starts at $-\infty$, hugs the lower
side of the negative real axis, then circles the origin in the
positive (counter-clockwise) direction before returning to
$-\infty$ along the upper side of the negative real axis. By
deforming the contour, one can readily check that~(\ref{qLaplace})
and~(\ref{qHankel}) are equivalent if $\Re(\beta)<0$ and $\beta$
is not a negative integer.  It follows~\cite[Cor.~2.5]{BraDi} that
if $\alpha_j$ is real for each $j=0,1,\dots,m$ and $\beta<0$, then
$q^*$ is positive, convex, and log-convex.  If $\beta$ is a
non-negative integer, then we can close the contour and deform it
into the unit circle, so that then $q^*$ is a polynomial of degree
$\beta$.  Otherwise, the singularity of the integrand is a branch
point at $z=0$, and $q^*$ is a transcendental function.  A
detailed treatment of the various properties of $q^*$, including
asymptotic expansions, behavior at $0$, and a representation of
$q^*(u)$ as an exponential of a Hellinger type
integro-differential operator acting on the monomial $u^\beta$ is
given in~\cite{BraDi}.

The main interest in $q^*$ stems from the fact that it provides a
$C^\infty$ solution to the difference differential
equation~(\ref{+ddereal}).  More specifically, combining
Proposition~2.2 and Corollary~2.4 of~\cite{BraDi} yields the
following result.
\begin{Prop}\label{prop:q*works}
The function $q^*$ of Definition~\ref{def:q*} satisfies the
difference differential equation~\textup{(\ref{+ddereal})} and the
asymptotic formula
\begin{equation}
   q^*(u)\sim u^\beta, \qquad 0<u\to\infty.
\end{equation}
\end{Prop}
We recall~\cite[Theorem 3]{BraDi} that there is at most one
function $q$ satisfying the difference differential
equation~(\ref{+ddereal}) and the additional condition that $q(u)$
is asymptotic to a fixed power of $u$ for large real values of the
argument $u$. Combining this result with
Proposition~\ref{prop:q*works} yields the following uniqueness
result, which we record for future use:
\begin{Thm}\label{thm:asym} Suppose that $q$ satisfies the
difference differential equation~\textup{(\ref{+ddereal})} and the
additional requirement that there exists a complex number $\tau$
such that
\[
   q(u) \sim u^\tau,\qquad 0<u\to\infty.
\]
Then $\tau=\beta$ and $q=q^*$, where $q^*$ is as in
Definition~\ref{def:q*}.
\end{Thm}
In other words, the function $q^*$ of Definition~\ref{def:q*} is
the unique solution to the difference differential
equation~(\ref{+ddereal}) in the class of functions asymptotic to
a monomial at infinity. Theorem 4 of~\cite{BraDi} strengthens
Theorem~\ref{thm:asym} by widening the class of admissible
functions to those functions which are polynomially bounded (i.e.\
majorized by a polynomial) in the right half-plane. The proof
in~\cite{BraDi} of this latter result is unnecessarily
complicated, applying residue theory to the inverse Laplace
transform of a suitably high order derivative of a supposed
solution.  Before proceeding further, it may be of interest to
give a simpler proof which has the additional advantage of
requiring only that the supposed solution be majorized by a
polynomial on the positive real half-line.

\begin{Thm}\label{thm:poly}
Suppose that $q$ satisfies the difference differential
equation~\textup{(\ref{+ddereal})}.  Suppose further that there
exists a real number $r$ such that
\begin{equation}
   q(u) = \bigoh(u^r),\qquad 0<u\to\infty.
   \label{polygrowth}
\end{equation}
Then there exists a complex number $A$ such that for all $u>0$,
$q(u)=Aq^*(u)$, where $q^*$ is given by Definition~\ref{def:q*}.
The multiplicative constant $A$ may of course depend on the
parameters $\alpha_j$, $v_j$ $(j=0,1,2,\dots,m)$ but not on $u$.
\end{Thm}

\demo{Proof}  For any complex number $A$, Theorem~\ref{thm:asym}
implies that $q=Aq^*$ is the unique function satisfying both the
difference differential equation~(\ref{+ddereal}) and the
asymptotic formula
\begin{equation}
   q(u) \sim Au^{\beta},\qquad 0<u\to\infty.\label{qinf}
\end{equation}
Accordingly, it is enough to show that our hypotheses imply the
asymptotic formula~(\ref{qinf}) holds for some complex number $A$,
or equivalently, that $$
   \lim_{u\to\infty} u^{-\beta}q(u)
$$
exists.

We can write the difference differential equation~(\ref{+ddereal})
in the form
\[
   \left(u^{-\beta}q(u)\right)'
   = u^{-\beta-1}\sum_{j=1}^m \alpha_j\{q(u+v_j)-q(u)\},
   \qquad u>0.
\]
We have dropped the term in the sum that corresponds to $j=0$.  It
vanishes because $v_0=0$.  It follows that for $u>0$,
\begin{align}
   u^{-\beta}q(u) &= q(1) + \int_1^u t^{-\beta-1}\sum_{j=1}^m \alpha_j
   \{q(t+v_j)-q(t)\}\,dt\nonumber\\
   &= q(1) + \sum_{j=1}^m \alpha_j\int_1^u t^{-\beta-1}
   \int_0^{v_j} q'(x+t)\,dx\,dt.
   \label{qdoubleint}
\end{align}
Let $c:=\Re(\beta)$, and suppose that the polynomial growth
requirement~(\ref{polygrowth}) holds. Since increasing $r$ only
weakens the hypothesis, there is no harm in assuming that $r=
n+c+1/2$ for some non-negative integer $n$. The stipulation that
$r-c$ lies halfway between two consecutive integers is not
strictly necessary; however, as it precludes the possibility that
an integer power of $u$ multiplied by $u^{r-c}$ will integrate to
a logarithm as opposed to a power of $u$, the number of cases to
be considered in our subsequent estimates is conveniently reduced.

By the difference differential equation~(\ref{+ddereal}) and our
hypothesis, we have
\[
   |uq'(u)| \le \sum_{j=0}^m |\alpha_j||q(u+v_j)| = \bigoh(u^r),
   \qquad 0<u\to\infty,
\]
and hence $q'(u) = \bigoh(u^{r-1}).$  Using this bound in the
integral~(\ref{qdoubleint}), we conclude that
\[
   u^{-\beta}q(u) = q(1) + \bigoh\left(\int_1^u t^{-c-1} t^{r-1}\,dt \right)
    = \bigoh(1) + \bigoh(u^{r-c-1}),
\]
and hence $q(u) = \bigoh(u^c)+\bigoh(u^{r-1})$. By repeating the
previous steps if necessary, we arrive at the growth estimate
\[
   q(u) = \bigoh(u^c),\qquad 0<u\to\infty.
\]
Using this latter bound in the difference differential
equation~(\ref{+ddereal}) yields
\[
   |uq'(u)| \le \sum_{j=0}^m |\alpha_j||q(u+v_j)|
    = \bigoh(u^c),\qquad 0<u\to\infty,
\]
and hence $q'(u) = \bigoh(u^{c-1}).$  It follows that the
integrals
\[
   \int_1^\infty t^{-\beta-1}\int_0^{v_j}q'(x+t)\,dx\,dt
   \qquad
   (j=0,1,2,\dots,m)
\]
all converge.  Therefore, we may let $u\to+\infty$
in~(\ref{qdoubleint}), from which we infer that the limit
\[
   A:= \lim_{u\to\infty} u^{-\beta}q(u) = q(1) +
   \sum_{j=1}^m\alpha_j\int_1^\infty t^{-\beta-1}
   \int_0^{v_j}q'(x+t)\,dx\,dt
\]
exists. In other words,
\[
   q(u) \sim A u^{\beta},\qquad 0<u\to\infty,
\]
and by our initial remarks, the proof is complete. \qed
\enddemo

In \S\ref{sect:exp}, we show that the conclusion of
Theorem~\ref{thm:poly} holds with~(\ref{polygrowth}) replaced by
the weaker hypothesis that $q$ be majorized by some fixed power of
an exponential function.  Before proving this stronger result in
full generality, it is instructive to examine an important special
case which is used in the proof of the general case.

\section{The Equation $ (uq(u))' = \kappa q(u)-\kappa q(u+1)$}
\label{sect:sieveaux}

Let $\kappa$ be a positive real number.  The difference
differential equation
\begin{equation}
   (uq(u))' = \kappa q(u)-\kappa q(u+1), \qquad u>0,
   \label{qkappa}
\end{equation}
was introduced by Iwaniec~\cite{Iw} as an adjoint equation for the
upper and lower bounding sieve functions of Rosser's sieve, and is
a special case of~(\ref{+ddereal}) with $m=1$, $v_1=1$,
$\alpha_0=\kappa-1$, $\alpha_1=-\kappa$ and
$\beta=\alpha_0+\alpha_1=-1.$  The role of $\kappa$ in sieve
theory is to measure the average number of residue classes being
deleted for each prime used in the sifting.
Equation~(\ref{qkappa}) and its close cousin with the minus sign
replaced by a plus sign have been studied fairly extensively in
the context of sieve
theory~\cite{DBrad0,DBrad1,DHR1II,DHR2,DHR3,DHR4,DHR5,DHR6,Iw}.
The case $\kappa=1$ of~(\ref{qkappa}) also occurs in the problem
of determining the asymptotic average size of the largest prime
factor of a random integer and the longest cycle of a random
permutation~\cite{Gol,Gonch,KTP,SL}.  More specifically, the
average cycle length of the longest cycle in a permutation on $n$
symbols~\cite{SL} and the average number of digits in the largest
prime factor of an $n$-digit number~\cite{KTP} are both asymptotic
to $nq_1'(1)/q_1(1)$, where
\[
   q_1(u):=\int_0^\infty \exp\left\{-ux-
   \int_0^x \frac{1-e^{-t}}{t}\,dt\right\}\,dx,\qquad u>0,
\]
satisfies~(\ref{qkappa}) with $\kappa=1$.  We also
have~\cite{Wh1,Wh2}
\[
   \#\big\{n\in\Z:1\le n\le x, P_2(n)\le P_1(n)^{1/u}\big\}
   = xe^{\gamma}q_1(u)+\bigoh\bigg(\frac{x}{\log x}\bigg),\qquad x\to\infty,
\]
where $P_j(n)$ denotes the $j$th largest prime factor of $n$ and
$\gamma$ is Euler's constant.

By the results of \S\ref{sect:poly}, we know that the polynomially
bounded solutions to~(\ref{qkappa}) are all constant multiples of
$q^*$, where here
\begin{equation}
   q^*(u)=\il\exp\bigg\{-ux-\kappa\int_0^x\frac{1-e^{-t}}{t}\,dt
   \bigg\}\,dx,\qquad u>0.
   \label{qsoln}
\end{equation}
It turns out that any solution of~(\ref{qkappa}) which is not a
constant multiple of~(\ref{qsoln}) must oscillate wildly. More
precisely, we show that if $q$ satisfies~(\ref{qkappa}) with
$\kappa>0$, and if there is no constant $A$ for which $q=Aq^*$,
then $q(u)\ne \bigoh(\exp(\lambda u))$ for any fixed $\lambda>0$,
and for all $u>0$, $q$ changes sign infinitely often in
$[u,\infty)$.  To prove this, we first establish the following
special case of our main uniqueness result.
\begin{Lem}\label{lem:beta=-1}
Suppose that $q$ satisfies the difference differential
equation~\textup{(\ref{+ddereal})} with $\beta=-1$.  Suppose
further that there exists a positive real number $\lambda$ such
that
\[
   q(u) = \bigoh(e^{\lambda u}),\qquad 0<u\to\infty.
\]
Then there exists a complex number $A$ such that for all $u>0$,
$q(u)=Aq^*(u)$, where $q^*$ is given by Definition~\ref{def:q*}
with $\beta=-1$, i.e.
\[
   q(u) = A\il \exp\bigg\{-ux+\sum_{j=1}^m\alpha_j
   \int_0^x \frac{1-e^{-v_j t }}{t}\,dt\bigg\}\,dx,\qquad u>0.
\]
\end{Lem}

\demo{Proof} Partial summation enables us to rewrite the
differential difference equation~(\ref{+ddereal}) in the form
\[
   (uq(u))' = \sum_{j=1}^m c_j \{q(u+v_j)-q(u+v_{j-1})\},\qquad u>0,
\]
where
\[
   c_j = \sum_{i=j}^m \alpha_i,\qquad j=1,2,\dots,m.
\]
Now integrate.  There exists a complex number $A$ such that for
all $u>0$,
\begin{equation}
   uq(u) = A+\sum_{j=1}^m c_j\int_{v_{j-1}}^{v_j} q(u+t)\,dt.
\label{uqu}
\end{equation}
Let
\[
   M := \frac{4}{\lambda}\sum_{j=1}^m |c_j| e^{\lambda v_j}.
\]
The exponential growth requirement implies that there exist
$u_0\ge M$ and $B>0$ such that for all $u\ge u_0$, $|q(u)|\le
Be^{\lambda u}.$ Inserting this inequality into~(\ref{uqu}) yields
\[
   |uq(u)|\le |A|+ B\sum_{j=1}^m |c_j|\int_{v_{j-1}}^{v_j}
   e^{\lambda(u+t)}\,dt,   \qquad u\ge u_0.
\]
It follows that
\begin{equation}
   |q(u)|\le \frac{|A|}{u} + \frac{Be^{\lambda u}}{\lambda u}\sum_{j=1}^m |c_j|e^{\lambda v_j}
   \le \frac{|A|}{u} + \frac{Be^{\lambda u}}{4},
   \qquad u\ge u_0.
   \label{qbound}
\end{equation}
Now we claim that $|q(x)|\le 2 |A|/x$ for $x\ge u_0$.  To prove
the claim, fix $x\ge u_0$ and observe that if $Be^{\lambda x}/4\le
|A|/x$,
then~(\ref{qbound}) gives $|q(x)|\le 2|A|/x.$
On the other hand, if $Be^{\lambda x}/4 > |A|/x$,
then we must have $Be^{\lambda u}/4 > |A|/u$ for all $u\ge x$.
Since~(\ref{qbound}) holds for all $u\ge u_0$, we then get
\[
   |q(u)|\le \frac{Be^{\lambda u}}{2},\qquad u\ge x.
\]
We can now insert this latter inequality back into~(\ref{uqu}) and
repeat the previous reasoning with $B$ replaced by $B/2$.  In
general, if $n$ is a positive integer such that
\[
   \frac{Be^{\lambda x}}{4\cdot2^{n-1}}>\frac{|A|}{x},
\]
then iterating the previous argument $n$ times will yield
\begin{equation}
   |q(u)|\le \frac{|A|}{u} + \frac{Be^{\lambda u}}{4\cdot2^n},\qquad u\ge x.
   \label{almost-there}
\end{equation}
Note that the bound~\eqref{almost-there} is valid for $u\ge x$,
where $x$ is the same as above. If $A\ne 0$, let $n$ be the least
positive integer such that
\[
   \frac{Be^{\lambda x}}{4\cdot2^n}\le \frac{|A|}{x}.
\]
Then~(\ref{almost-there}) gives
\begin{equation}
   |q(x)|\le \frac{2|A|}{x}.
\label{qxbound}
\end{equation}
Since $x\ge u_0$ is arbitrary,~(\ref{qxbound}) must hold for all
$x\ge u_0$.  On the other hand, if $A=0$, then letting
$n\to\infty$ in~(\ref{almost-there}) shows that $q(u)=0$ for all
$u\ge x$. But again, since $x\ge u_0$ is arbitrary, we must have
$q(u)=0$ for all $u\ge u_0$.  In other words,~(\ref{qxbound}) also
holds for all $x\ge u_0$ when $A=0$.   This establishes the claim.
Invoking Theorem~\ref{thm:poly} now completes the proof. \qed
%
\enddemo

We can now show that any solution to~(\ref{qkappa}) with
$\kappa>0$ must be wildly oscillatory.
%
\begin{Thm} Suppose that $q$ is a solution to the difference differential
equation~\textup{(\ref{qkappa})} with $\kappa>0$.  Suppose further
that $q$ is not a constant multiple of $q^*$, where $q^*(u)$ is
given by~\textup{(\ref{qsoln})}. Then $q(u)\ne \bigoh(\exp(\lambda
u))$ for any fixed $\lambda >0$, and for all $u>0$, $q$ changes
sign infinitely often in $[u,\infty)$.
\end{Thm}

\demo{Proof} Note that $\beta=-1$ in~(\ref{qkappa}).  The
conclusion $q(u)\ne \bigoh(\exp(\lambda u))$ for all $\lambda>0$
is an immediate consequence of Lemma~\ref{lem:beta=-1}.
 In particular
there exists a strictly increasing sequence of positive real
numbers $u_1<u_2<\dots$ satisfying $\lim_{n\to\infty}u_n=\infty$
and $|q(u_n)| > \exp(u_n)$ for all positive integers $n$.
Integrating the difference differential equation~(\ref{qkappa})
yields
\begin{equation}
   uq(u) = A - \kappa\int_u^{u+1} q(t)\,dt,\qquad u>0,\label{4.4}
\end{equation}
where $A$ is a constant of integration which may depend on
$\kappa$, but not on $u$.  There is a positive integer $N$ such
that if $n>N$ then $u_n\ge 1$ and $\exp(u_n) > |A|$.  Let $n>N$
and set $u=u_n$. Suppose that $q(u)>0$. Then~(\ref{4.4}) yields
\[
   \kappa\int_u^{u+1} q(t)\,dt = A - uq(u) < |A| - e^u < 0.
\]
Since $\kappa>0$, it follows that
\[
   \int_u^{u+1} q(t)\,dt <0,
\]
and hence $q$ must change sign from positive to negative somewhere
in the half-open, half-closed interval $(u,u+1]$. If, on the other
hand, $q(u)<0$, then~(\ref{4.4}) yields
\[
   \kappa\int_u^{u+1} q(t)\,dt = A - uq(u) > -|A| + e^u > 0,
\]
and hence $q$ must change sign from negative to positive somewhere
in $(u,u+1]$.  Thus, we have shown that $q$ has a sign change
beyond $u_n$ for each $n>N$, and the proof is complete.\qed
\enddemo

\section{The Adjoint Relation}
\label{sect:adj}

We'd like to remove the restriction $\beta=-1$ in
Lemma~\ref{lem:beta=-1}.
A careful examination of the proof of Lemma~\ref{lem:beta=-1}
suggests that the iterative argument employed therein succeeded
because the lengths of the ranges of integration in the
representation~(\ref{uqu}) are all independent of $u$.  Note
that~\eqref{uqu} is equivalent to
\[
   uq(u) = A + \sum_{j=1}^m \alpha_j\int_{u-v_j}^u q(t+v_j)\,dt.
\]
In order to recover this property in the case when $\beta\ne-1$,
let us seek a function $p$ such that for some constant $A$, the
equation
\begin{equation}
   up(u)q(u) = A + \sum_{j=1}^m \alpha_j\int_{u-v_j}^u
   p(t)q(t+v_j)\,dt
   \label{adjeqn}
\end{equation}
holds for all sufficiently large values of $u$.  Let us suppose
further that beyond some point, $p$ is differentiable.
Differentiating~(\ref{adjeqn}) with respect to $u$ reveals that
\[
   uq'(u)p(u)+up'(u)q(u)+p(u)q(u)=p(u)\sum_{j=1}^m\alpha_jq(u+v_j)
   -q(u)\sum_{j=1}^m\alpha_jp(u-v_j).
\]
If we now substitute the right hand side of the difference
differential equation~(\ref{+ddereal}) for $uq'(u)$ in this latter
equation, the first sum on the right hand side can be cancelled
with the corresponding sum on the left. The result is
\[
   \alpha_0p(u)q(u)+up'(u)q(u)+p(u)q(u)=-q(u)\sum_{j=1}^m\alpha_jp(u-v_j),
\]
which is certainly the case if $p$ satisfies the delay
differential equation~(\ref{-ddereal}).
Conversely, it is clear that if $q$ satisfies~(\ref{+ddereal}) and
$p$ satisfies~(\ref{-ddereal}), then~(\ref{adjeqn}) must also hold
for all sufficiently large values of $u$.

\section{The Function $u\mapsto p(u,a,b)$}
\label{sect:pasym}

Our plan is to use the adjoint equation~(\ref{adjeqn}) to deduce
the asymptotic behavior of $q(u)$ for large $u$ using only a very
weak estimate on the growth rate of $q(u)$. In order to do this,
we need reasonably precise knowledge concerning the rate of growth
of a non-trivial solution $p(u)$ to the delay differential
equation~(\ref{-ddereal}) as $u$ increases without bound.  Most of
what we need can be found in~\cite{LF}, where  the existence of an
asymptotic expansion for $p$ for the case (in our notation)
$\Re(\alpha_0)<0$ is proved, and the first term of the expansion
is determined explicitly.  Nevertheless, it is worthwhile to
obtain some additional results, and for this it is useful to
extend Wheeler's~\cite{Wh1,Wh2} more comprehensive treatment of
the $m=2$ case to arbitrary $m$ and complex $\alpha_j$. Therefore,
with regard to the choice of boundary condition and the location
and classification of discontinuities, our development in this
section more closely parallels that of~\cite{Wh2}.  However, to
provide an explicit formula for all coefficients in the asymptotic
expansion, we find it more convenient to extend the technique used
in~\cite{LF}.  To facilitate comparison with our results, we set
$a=1+\alpha_0$ and let $b$ denote the vector
$(\alpha_1,\alpha_2,\dots,\alpha_m)$.

Let
\[
   C_0 = \prod_{j=1}^m \big(v_j e^{\gamma}\big)^{-\alpha_j},
\]
where, as usual,
\[
   \gamma = \lim_{n\to\infty} \bigg(\sum_{j=1}^n \frac{1}{j}-\log
   n\bigg)
\]
is Euler's constant.  Consider a particular solution
$p(u)=p(u,a,b)$ of the delay differential
equation~(\ref{-ddereal}), which is defined uniquely for all real
$u$ by the following six conditions:
\begin{subequations}
\begin{align}
   & (up(u))'=-\sum_{j=0}^m \alpha_j p(u-v_j);\label{delay}\\
   & p(u) =0 \mbox{ whenever $u\le 0$};\label{flat}\\
   & p(u) = \frac{C_0 }{\G(1-a)} u^{-a}\mbox{ if $0<u\le
   v_1$};
            \label{initcond}\\
   & p(u) \mbox{ is continuous for $u>0$ when
           $\Re(a)<1$};\label{cts}\\
   & p(u,a+1,b)
   =\frac{d}{du}p(u,a,b);\label{shift}\\
   & p \mbox{ is continuous from the left at every point}\label{leftcts}.
\end{align}
\end{subequations}
We need to prove that a unique such function exists.  First assume
$\Re(a)<1$.  The function defined by~(\ref{initcond}) is
integrable for $0< u\le v_1$ and the  delay differential
equation~(\ref{delay}) can be rewritten as
\[
   \left(u^a p(u)\right)' = -u^{a-1}\sum_{j=1}^m
   \alpha_j p(u-v_j),
\]
which can be integrated forward on successive intervals.  That is,
for each positive integer $n$, if $n<u\le n+1$ then
\begin{equation}
   p(u v_1) = u^{-a}\bigg\{n^a p(nv_1)
        -  \int_n^u t^{a-1} \sum_{j=1}^m \alpha_j
        p(tv_1-v_j)\,dt\bigg\}.
\label{pintegrated}
\end{equation}
This establishes uniqueness and every condition
except~(\ref{shift}).  To see that~(\ref{shift}) holds, simply
differentiate~(\ref{delay}) and~(\ref{initcond}) and apply the
functional equation for the gamma function.  For the case
$\Re(a)\ge 1$, we use the fact that if $n$ is a positive integer,
then
\[
   p(u,a,b)=\bigg(\frac{d}{du}\bigg)^n
   p(u,a-n,b),
\]
which follows from~(\ref{shift}).

From~(\ref{pintegrated}), it follows that $p$ has a two-sided
derivative at every point $x$ for which $p$ is continuous at each
of the points $x-v_j$, $1\le j\le m$.  Otherwise, we need to take
derivatives from the left in~(\ref{delay}) and~(\ref{shift}). This
we can do, again by~(\ref{pintegrated}), since $p$ is continuous
from the left at every point.

\begin{Thm}\label{thm:pDiscnties}
\noindent

\textup{(i)} If $\Re(a)$ is not an integer, then the set of points
at which $p$ is discontinuous is precisely $\{nv_j: n\in\Z,\; 0\le
n<\Re(a),\; 1\le j\le m\}$.  At these points, $p$ has a finite
limit from the left, and is unbounded from the right.  Moreover,
if $n$ is a non-negative integer, $\Re(a)>n$, and $1\le j\le m$,
then
\[
   p(u)\sim
   \frac{(-\alpha_j)^n C_0}{n!\, v_j^n\,\G(n-a+1)}
   (u-nv_j)^{n-a},\qquad u\to nv_j+,
\]
whereas $\ds\lim_{u\to nv_j-} p(u)=p(nv_j)$ is finite.

\textup{(ii)} If $\Re(a)$ is a non-negative integer, but $a$ is
not itself an integer, then the set of points at which $p$ is
discontinuous is precisely $\{nv_j : n\in\Z,\; 0\le n\le \Re(a),\;
1\le j\le m\}$.  The discontinuities consist of two types.  If $n$
is a non-negative integer, $\Re(a)>n$, and $1\le j\le m$, then as
in \textup{(i)},
\[
   p(u)\sim
   \frac{(-\alpha_j)^n C_0}{n!\, v_j^n\,\G(n-a+1)}
   (u-nv_j)^{n-a},\qquad u\to nv_j+,
\]
whereas $\ds\lim_{u\to nv_j-} p(u)=p(nv_j)$ is finite. If
$n=\Re(a)$, then $p$ is bounded in a neighbourhood of $nv_j$, but
$\ds\lim_{u\to nv_j+} p(u)$ does not exist.

\textup{(iii)} If $a$ is a positive integer, then $p$ has only
finite jump discontinuities.  These can only occur at the points
$nv_j$, where $n$ is an integer, $1\le n\le a$, and $1\le j\le m$.

\textup{(iv)} Finally, $p(u,0,b)$ has a finite jump discontinuity
at $u=0$ and no other discontinuities.
\end{Thm}

\demo{Sketch of Proof} The statements of the Theorem are true if
$\Re(a)<1$ as we already noted.  The other cases can be
established by induction, using the equation
\begin{equation}
   p(u,a+1,b) = \frac{d}{du}p(u,a,b)
   = -\frac{a}{u}p(u,a,b)-\frac{1}{u}\sum_{j=1}^m\alpha_j
   p(u-v_j,a,b),
\label{pdiscont}
\end{equation}
which is certainly valid if $u$ is not an integer multiple of
$v_j$, $1\le j\le m$. For details, see~\cite[Theorem 1]{Wh2}\qed
\enddemo

\begin{Rem} The situation in (iii) is somewhat unsatisfactory,
as in general all we can say is that if $a$ is a positive integer
and $u$ is a point of discontinuity, then necessarily $u=nv_j$ for
some positive integer $n\le a$ and $1\le j\le m$. Although it is
typically the case that all such points are indeed points of
discontinuity---for example, this is indeed true if the $m$
positive real numbers $v_1,\dots,v_m$ are linearly independent
over the integers---it may happen that for certain values of the
parameters $\alpha_j$, $v_j$, the jumps arising from different
terms in~(\ref{pdiscont}) cancel each other.  As an example,
consider the case $m=2$, $a=2$, $b=(\alpha_1,\alpha_2)=(1,-2)$,
$v_1=1$, $v_2=2$.  Then $C_0 = 4 e^\gamma$.  We find that the jump
at $2v_1+$ is exactly cancelled by the jump at $v_2+$, so that
$p(2+,2,b)=e^\gamma = p(2-,2,b)$. As noted by Wheeler~\cite{Wh2},
such cancellations cannot occur when $m=1$.
\end{Rem}

We now turn to the problem of determining the behavior of $p(u)$
for large positive real $u$.  We take the following result of
Levin and Fainleib~\cite{LF} as our point of departure, recasting
it in our notation, and with our boundary
condition~(\ref{initcond}) etc.

\begin{Lem}[\cite{LF}, Lemma 1.3.1]\label{Lem:LF} Let $\Re(a)<1$,
and let $\varphi(u)=(u\log u)/v_m + \bigoh(u)$, $0<u\to \infty$.
We have the following behavior for $p(u)$ as $0<u\to\infty$.

\textup{(i)} If $\beta$ is a non-negative integer, then
$p(u)=\bigoh(\exp(-\varphi(u)))$.

\textup{(ii)} If $\beta$ is a negative integer, there exists a
polynomial $r_\beta$ of degree $-\beta-1$ such that $p(u) =
r_\beta(u)+\bigoh(\exp(-\varphi(u)))$.

\textup{(iii)} For any complex $\beta$ not an integer, there is an
asymptotic expansion $p(u)\sim \sum_{n\ge 0} c_n u^{-\beta-1-n}$.
\end{Lem}

Levin and Fainleib~\cite{LF} did not identify the polynomial
$r_\beta$ in (ii).  In addition, they gave only the first term in
the asymptotic expansion (iii). We shall determine all the
coefficients $c_n$ and the polynomial $r_\beta$ explicitly. We
shall also remove the restriction on $\Re(a)$.  To carry this out,
we need the following formula for the Laplace transform of $p$.

\begin{Thm}\label{thm:pLaplace} Let $\Re(a)<1$ and $\Re(s)>0$.  Then
\[
   \int_0^\infty e^{-su} p(u)\,du = s^\beta
   \exp\bigg\{-\sum_{j=1}^m \alpha_j \int_0^s
   \frac{1-e^{-tv_j}}{t}\,dt\bigg\}.
\]
If $\beta$ is a non-negative integer, the formula is valid for all
complex numbers $s$.
\end{Thm}

Theorem~\ref{thm:pLaplace} is actually just a restatement in our
notation of a corresponding result derived in~\cite{LF}.  Although
the region of convergence was not discussed there, both the
existence of the Laplace transform and its region of convergence
follow easily from the growth behavior of $p$ as given in
Lemma~\ref{Lem:LF}.


Let $n$ be a non-negative integer.  Recall the
polynomial~\cite[equation~(3.2)]{BraDi}
\[
   Q_n(u,b) = \bigg(\frac{\partial}{\partial
   z}\bigg)^n\bigg|_{z=0}
   \exp\bigg\{uz-\sum_{j=1}^m\alpha_j\int_0^z\frac{e^{v_jt}-1}{t}\,dt
   \bigg\}.
\]
Remembering that $b$ denotes the vector of coefficients
$(\alpha_1,\dots,\alpha_m)$, we let $Q_n(u,-b)$ be the polynomial
obtained by negating each $\alpha_j$ $(1\le j\le m)$ in
$Q_n(u,b)$.  With this notation, we can now restate
Lemma~\ref{Lem:LF} with the coefficients $c_n$ and the polynomial
$r_\beta$ explicitly identified.

\begin{Thm}\label{thm:pasym} Let $\varphi(u)$ be as in the statement of
Lemma~\ref{Lem:LF}.  As $0<u\to \infty$, we have the following
behavior for $p(u)$.

\textup{(i)} If $\beta$ is a non-negative integer, then
$p(u)=\bigoh(\exp(-\varphi(u)))$.

\textup{(ii)} If $\beta$ is a negative integer, then
\[
   p(u) = \sum_{n=0}^{-\beta-1}\frac{(-1)^n}{n!}
   \frac{ Q_n(0,-b)}{(-\beta-1-n)!}u^{-\beta-1-n}
   +\bigoh\big(e^{-\varphi(u)}\big).
\]

\textup{(iii)} For any complex number $\beta$ not an integer, we
have the asymptotic expansion
\[
   p(u)\sim \sum_{n=0}^\infty \frac{(-1)^n}{n!}
   \frac{Q_n(0,-b)}{\G(-\beta-n)}u^{-\beta-1-n}.
\]
\end{Thm}

\demo{Proof} In view of Lemma~\ref{Lem:LF}, there is nothing to
prove for (i), and the error term need not concern us in (ii). Fix
a positive integer $h$. It is sufficient to consider the case
$\Re(a)<1$ and $\Re(\beta)<1-h$, since as noted in~\cite{LF}, the
remaining cases reduce to the problem of differentiating
asymptotic expansions, and by~(\ref{delay}) this is legitimate.
Let $0<\e<1$ be such that $-\Re(\beta)-h-\e>-1$. Substituting
\[
   p(u) = \sum_{n=0}^{h-1} c_n u^{-\beta-1-n}
   +\bigoh(u^{-\Re(\beta)-h-\e}),\qquad 0<u\to\infty
\]
into the Laplace transform for $p$, we find that
\[
   \int_0^\infty e^{-su} p(u)\,du
    = \sum_{n=0}^{h-1}c_n\G(-\beta-n)s^{\beta+n}
   +\bigoh(s^{\Re(\beta)+h+\e-1}),\qquad s\to 0+.
\]
On the other hand, Theorem~\ref{thm:pLaplace} and the definition
of the polynomials $Q_n(0,-b)$ imply that
\[
   \int_0^\infty e^{-su} p(u)\,du\sim \sum_{n=0}^\infty
   \frac{(-1)^n}{n!}Q_n(0,-b)s^{\beta+n},
   \qquad s\to 0+.
\]
Comparing the two expressions for the Laplace transform, we infer
that
\[
   c_n=\frac{(-1)^n}{n!}\frac{Q_n(0,-b)}{\G(-\beta-n)},
\]
which proves (iii).  If $\beta$ is an integer and $\beta\ge -n$,
we see that $c_n=0$.   It follows that if $\beta$ is a negative
integer, then
\[
   r_\beta(u) = \sum_{n=0}^{-\beta-1}\frac{(-1)^n}{n!}
   \frac{Q_n(0,-b)}{\G(-\beta-n)}u^{-\beta-1-n},
\]
as stated in (ii). \qed
\enddemo

\section{Exponentially Bounded Solutions}
\label{sect:exp}

In this section we prove that, up to a multiplicative constant,
the difference differential equation~(\ref{+ddereal}) has a unique
solution in the class of functions majorized by a function of
exponential growth.  We first establish an inequality for certain
exponential integrals that arise in the proof.

\begin{Lem}\label{lem:IncompleteGamma}
Let $u$, $r$ and $\lambda$ be real numbers satisfying $\lambda
>r^+=\max(r,0)$ and $u>1$.  Then
\[
   \int_1^u e^{\lambda t}t^{-r}\,dt \le \frac{e^{\lambda
   u}u^{-r}}{\lambda -r^+}.
\]
\end{Lem}

\demo{Proof} Denote the integral by $I$.  An easy integration by
parts shows that
\begin{equation}
   I\le \frac{1}{\lambda}e^{\lambda u}u^{-r}+\frac{r}{\lambda}
   \int_1^u e^{\lambda t}t^{-r-1}\,dt.
\label{IP}
\end{equation}
If $r<0$, then the coefficient of the integral on the right hand
side of~(\ref{IP}) is negative, and so the corresponding term can
be dropped, yielding
\[
   I\le \frac{1}{\lambda}e^{\lambda u}u^{-r}.
\]
Since $r<0$ implies $r^+=0$, the claim is established in this
case.  On the other hand, if $r\ge 0$, then $t^{-r}\ge t^{-r-1}$
for $1\le t\le u$ and hence
\[
   I\le \frac{1}{\lambda}e^{\lambda u}u^{-r}+\frac{r}{\lambda}
   \int_1^u e^{\lambda t}t^{-r}\,dt
   = \frac{1}{\lambda}e^{\lambda u}u^{-r}+\frac{r}{\lambda}I.
\]
Since $\lambda>r\ge 0$, it follows that
\[
   I\le \frac{e^{\lambda u}u^{-r}}{\lambda -r}.
\]
Since $r\ge 0$ implies $r^+=r$, the claim is established in this
case also. \qed
\enddemo

\begin{Thm}\label{thm:exp}
Suppose that $q$ satisfies the difference differential
equation~\textup{(\ref{+ddereal})}.  Suppose further that there
exists a positive real number $\lambda$ such that
\[
   q(u) = \bigoh(e^{\lambda u}),\qquad 0<u\to\infty.
\]
Then there exists a complex number $A$ such that for all $u>0$,
$q(u)=Aq^*(u)$, where $q^*$ is given by Definition~\ref{def:q*}.
The multiplicative constant $A$ may of course depend on the
parameters $\alpha_j$, $v_j$ $(j=0,1,2,\dots,m)$ but not on $u$.
\end{Thm}

\demo{Proof} First, suppose that $\beta=n$ is a non-negative
integer.  The difference differential equation~(\ref{+ddereal})
implies that $q'(u)$ satisfies~\eqref{+ddereal} with $\alpha_0$
replaced by $\alpha_0-1$, and that
$q'(u)=\bigoh\left(u^{-1}e^{\lambda u}\right)$, $0<u\to\infty$.
Hence, we deduce that the $n+1^{st}$ derivative of $q$ satisfies
the conditions of Lemma~\ref{lem:beta=-1}. Accordingly, by
Lemma~\ref{lem:beta=-1} and Proposition~\ref{prop:q*works}, there
exists a complex constant $A$ such that
\[
   q^{(n+1)}(u)\sim A u^{-1},\qquad 0<u\to\infty.
\]
It follows that
\[
   q^{(n)}(u)\sim A\log u,\qquad 0<u\to\infty.
\]
Also, $q^{(n)}$ satisfies the difference differential
equation~(\ref{+ddereal}) with $\beta=0$.  Theorem~\ref{thm:poly}
now implies that $q^{(n)}$ is a constant multiple of $q^*$, where
$q^*$ is given by Definition~\ref{def:q*} with $\beta=0$.  It
follows that $A=0$ and $q^{(n)}$ is constant.  Hence, there exist
complex numbers $a_0,a_1,\dots,a_n$ such that for all $u>0$,
\[
   q(u) = \sum_{j=0}^n a_j u^j.
\]
More specifically, Theorem~\ref{thm:asym} implies that $q=a_nQ_n$,
where
\[
   Q_n(u) := \bigg(\frac{\partial}{\partial z}\bigg)^n\bigg|_{z=0}
   \exp\bigg\{uz-\sum_{j=0}^m \alpha_j\int_0^z\frac{e^{v_j t}-1}{t}\,dt\bigg\}
\]
is the polynomial obtained by setting $\beta=n$ in~(\ref{qHankel})
of Definition~\ref{def:q*} and deforming the contour into the unit
circle.  See~\cite{BraDi} for additional properties of the
polynomial $Q_n$.

For the general case, we assume that $\beta$ is a complex number,
but not a non-negative integer. Let $c=\Re(\beta)$. Since
increasing $\lambda$ only weakens the hypothesis, we may assume
that $K:=\lambda-\max(c+1,0)>0$. Let
\[
   M := \frac{16}{K}\sum_{j=1}^m |\alpha_j|e^{\lambda v_j}.
\]
By hypothesis, there exist $u_1\ge M$ and $B>0$ such that for all
$u\ge u_1$, $|q(u)|\le Be^{\lambda u}$. Let
$p(u)=\G(-\beta)p(u,a,b)$, where $p(u,a,b)$ is as
in~\S\ref{sect:pasym}.  Since $Q_0(0,-b)=1$,
Theorem~\ref{thm:pasym} implies that $p(u)\sim u^{-\beta-1}$ as
$0<u\to\infty$.  Therefore, there exists a positive real number
$u_2$ such that for all $u\ge u_2$, $\tfrac12 u^{-c-1}\le
|p(u)|\le 2u^{-c-1}$. There must also exist a complex number $A$
and a positive real number $u_3$ such that the adjoint relation
(cf.~\ref{adjeqn})
\begin{equation}
   up(u)q(u) = A + \sum_{j=1}^m \alpha_j\int_{u-v_j}^u
   p(t)q(t+v_j)\,dt
\label{innerprod}
\end{equation}
holds for all $u\ge u_3$.  Let $u_4 :=
\max(u_1,u_2+v_m,u_3,1+v_m)$. Then for all $u\ge u_4$, by
Lemma~\ref{lem:IncompleteGamma} we have
\begin{align*}
   |up(u)q(u)|&\le |A|+2B\sum_{j=1}^m |\alpha_j|
   \int_{u-v_j}^u t^{-c-1}e^{\lambda (t+v_j)}\,dt\\
      &\le |A| + 2 BK^{-1}u^{-c-1}\sum_{j=1}^m |\alpha_j|e^{\lambda
      (u+v_j)}\\
   &\le |A|+\frac{u_1}{8}B u^{-c-1}e^{\lambda u}.
\end{align*}
Let $A'=2 A$.  It follows that
\begin{equation}
   |q(u)|\le |A' u^{\beta}| + \frac{u_1}{4u}Be^{\lambda u}
   \le |A' u^{\beta}| + \tfrac14{B}e^{\lambda u},
   \qquad u\ge u_4.
\label{qll}
\end{equation}
Now let $u_0 := \max(u_4,c/\lambda)$, and fix $x\ge u_0.$  If
$\tfrac14 Be^{\lambda x} \le |A'x^{\beta}|$, then~(\ref{qll})
gives $|q(x)|\le 2|A' x^{\beta}|$.

On the other hand, if $\tfrac14B e^{\lambda x}>|A' x^{\beta}|$,
then we claim that $\tfrac14Be^{\lambda u}>|A' u^{\beta}|$ for all
$u\ge x$.  To see the claim, let $f(u):=\log B-\log 4+\lambda
u-\log |A'|-\Re(\beta)\log u.$  We have $f(x)>0$ and
$f'(u)=\lambda-c/u>0$ for all $u\ge x$.  Therefore, $f$ is
positive and strictly increasing on the interval $[x,\infty)$, and
this proves the claim.  Thus, in this case,~(\ref{qll}) gives
$|q(u)|\le \tfrac12Be^{\lambda u}$ for all $u\ge x$.  Now insert
this latter inequality back into the adjoint
relation~(\ref{innerprod}), i.e.~repeat the previous argument with
$B$ replaced by $B/2$.  We get
\begin{equation}
   |q(u)|\le |A' u^\beta| + \tfrac18Be^{\lambda u}, \qquad u\ge x.
\label{q18}
\end{equation}
In general, if $n$ is a positive integer such that
\[
   \frac{Be^{\lambda x}}{4\cdot 2^{n-1}} > |A' x^\beta|,
\]
then $n$ iterations of the preceding argument will yield
\begin{equation}
   |q(u)|\le |A' u^\beta|+\frac{Be^{\lambda u}}{4\cdot 2^n},
   \qquad u\ge x.
\label{qpow2}
\end{equation}
Therefore, if $A' \ne 0$, let $n$ be the least positive integer
such that
\[
   \frac{Be^{\lambda x}}{4\cdot 2^n}\le |A' x^\beta|.
\]
Then~(\ref{qpow2}) implies that
\begin{equation}
   |q(x)|\le 2|A' x^\beta|.
\label{qx}
\end{equation}
Since $x\ge u_0$ is arbitrary, the bound~(\ref{qx}) must hold for
all such $x$.  On the other hand, if $A'=0$, then letting
$n\to\infty$ in~(\ref{qpow2}) shows that $q(u)=0$ for all $u\ge
x$.  But again, $x\ge u_0$ is arbitrary, so we must have $q(u)=0$
for all $u\ge u_0$.  In other words, the bound~(\ref{qx}) holds
for all $x\ge u_0$ and all complex numbers $A'$.  Invoking
Theorem~\ref{thm:poly} completes the proof. \qed
\enddemo

\section{Final Remarks}
\label{sect:final}  The introduction of an adjoint equation as an
aid to studying solutions of difference differential equations
goes back at least to de Bruijn's study of the Buchstab
function~\cite{dB1}. There and subsequently (see
eg.~\cite{HM,CG,Wh1,Wh2}), the adjoint equation is used to deduce
information about a solution to a special case of~(\ref{-ddereal})
from a solution to the corresponding special case
of~(\ref{+ddereal}).  In this paper, we turned the process around,
deducing information about a solution to~(\ref{+ddereal}) from the
adjoint equation and a solution to~(\ref{-ddereal}).

It is interesting to compare and contrast the behavior of $p$ and
$q$ near zero and infinity.  As in~\S\ref{sect:exp} let
$p(u)=\G(-\beta)p(u,a,b)$, where $p(u,a,b)$ is as
in~\S\ref{sect:pasym}.  Assume that $\beta$ is not a non-negative
integer.  Then Theorem~\ref{thm:pasym} implies that
\[
   p(u)\sim\sum_{n=0}^\infty (-1)^n \binom{-\beta-1}{n} Q_n(0,-b)
   u^{-\beta-1-n},\qquad 0<u\to\infty,
\]
whereas Theorem 5 of~\cite{BraDi} states that
\[
   q(u) \sim \sum_{n=0}^\infty \binom{\beta}{n} Q_n(0,b)
   u^{\beta-n},\qquad 0<u\to\infty.
\]
If in addition, $\Re(\alpha_0)<0$, then Theorem 7(ii)
of~\cite{BraDi} can be restated as
\[
    q(u)\sim u^{\alpha_0} \frac{\G(-\alpha_0)}{\G(-\beta)}
    \prod_{j=1}^m \big(v_je^\gamma\big)^{\alpha_j},
    \qquad u\to 0+,
\]
whereas~(\ref{initcond}) implies that
\[
   p(u) = u^{-\alpha_0-1}\frac{\G(-\beta)}{\G(-\alpha_0)}
   \prod_{j=1}^m \big(v_j
   e^\gamma\big)^{-\alpha_j},\qquad 0<u\le v_1.
\]
In particular,
$\lim_{u\to0+}up(u)q(u)=\lim_{u\to+\infty}up(u)q(u)=1$.

\section*{Acknowledgment}
It is a pleasure to thank the referee for a careful examination of
the paper, and for several insightful suggestions which led to
improvements in the exposition.

\end{document}